\documentclass[12pt]{article}
\usepackage{amsmath}
\usepackage{amssymb}
\usepackage{amscd}
\usepackage{latexsym}
\usepackage{amsthm}
\usepackage{epsfig}
\usepackage{amsfonts}
\usepackage{enumerate}
\usepackage{xypic}

\newtheorem{theorem}{Theorem}[section]
\newtheorem*{theorem*}{Theorem}
\newtheorem{lemma}[theorem]{Lemma}
\newtheorem{proposition}[theorem]{Proposition}
\newtheorem{corollary}[theorem]{Corollary}

\newtheorem{conjecture}[theorem]{Conjecture}

{
\theoremstyle{plain}

\newtheorem{example}[theorem]{Example}

\newtheorem{remark}[theorem]{Remark}
}

\newcommand{\stab}{\operatorname{stab}}

\renewcommand{\Im}{\operatorname{Im}}

\newcommand{\codim}{\operatorname{codim}}
\renewcommand{\dim}{\operatorname{dim}}

\newcommand{\Sym}{\operatorname{Sym}}
\newcommand{\rk}{\operatorname{rk}}
\newcommand{\crk}{\operatorname{crk}}

\newcommand{\C}{{\mathbb{C}}}
\newcommand{\Z}{{\mathbb{Z}}}

\newcommand{\R}{{\mathbb{R}}}

\renewcommand{\O}{\mathcal{O}}

\newcommand{\otn}{\{1,\ldots,n\}}

\newcommand{\bigmid}{\hs\Big{|}\hs}
\newcommand{\subs}{\subseteq}
\newcommand{\sups}{\supseteq}
\newcommand{\hookto}{{\hookrightarrow}}
\renewcommand{\iff}{\Leftrightarrow}

\newcommand{\becircled}{\mathaccent "7017}
\newcommand{\hs}{\hspace{3pt}}

\newcommand{\ith}{i^{\text{th}}}

\newcommand{\D}{\Delta}
\renewcommand{\i}{\iota}
\renewcommand{\a}{\alpha}
\renewcommand{\b}{\beta}

\renewcommand{\cot}{T^*\!\Cn}
\newcommand{\Cnd}{(\Cn)^*}
\newcommand{\Cn}{\mathbb A^{\! n}}

\newcommand{\M}{\mathfrak{M}}
\newcommand{\Ma}{\M}
\newcommand{\X}{\mathfrak{X}}

\newcommand{\muc}{\mu}

\newcommand{\phc}{\Phi}

\newcommand{\mil}{\mu^{-1}(\lambda)}

\newcommand{\mcinv}{\muc^{\! -1}(0)}
\newcommand{\minv}{\mu^{\! -1}(0)}
\newcommand{\mcinvast}{\muc^{\! -1}(0)^{\a-\text{st}}}

\newcommand{\pinv}{\pi^{-1}}

\newcommand{\bigmod}{{\Big/\!\!\!\!\Big/}}
\renewcommand{\mod}{{/\!\!/}}
\newcommand{\mmod}{{\!/\!\!/\!\!/\!\!/}}

\newcommand{\Tdc}{T^d_{\C}}
\newcommand{\Tkc}{T^k_{\C}}

\newcommand{\Thc}{\hat T_\C}

\newcommand{\Tk}{T^k}
\newcommand{\Tn}{T^n}
\newcommand{\Td}{T^d}

\newcommand{\Trkf}{T^{\rk F}}
\newcommand{\Tcrkf}{T^{\crk F}}

\newcommand{\tk}{\mathfrak{t}^k}
\newcommand{\tn}{\mathfrak{t}^n}

\newcommand{\td}{\mathfrak{t}^d}
\newcommand{\tf}{\mathfrak{t}^{F}}
\newcommand{\tfc}{\mathfrak{t}^{F^c}}

\newcommand{\tkd}{(\tk)^*}

\newcommand{\tnd}{(\tn)^*}
\newcommand{\tdd}{(\td)^*}

\newcommand{\TG}{T^G}
\newcommand{\TF}{T^F}

\newcommand{\A}{\mathcal{A}}

\newcommand{\At}{\tilde\A}
\newcommand{\Aft}{\tilde\Af}
\newcommand{\Af}{\mathcal{A}_F}
\newcommand{\AF}{\mathcal{A}^F}
\newcommand{\MAft}{\M(\Aft)}
\newcommand{\MAt}{\M(\At)}
\newcommand{\MAtl}{\M_{\lambda}(\At)}
\newcommand{\XA}{\X(\A)}
\newcommand{\MA}{\M(\A)}
\newcommand{\MAF}{\M(\AF)}
\newcommand{\MAf}{\M(\Af)}
\newcommand{\LAft}{\mathcal{L}(\Aft)}
\newcommand{\LAt}{\mathcal{L}(\At)}
\newcommand{\MAopen}{\becircled\MA}
\newcommand{\MAFopen}{\becircled\MAF}

\newcommand{\Fpbar}{\bar{\mathbb F}_{\! p}}
\newcommand{\Fqbar}{\bar{\mathbb F}_{\! q}}
\newcommand{\da}{\Delta_\A}
\newcommand{\bcsda}{\operatorname{bc}_\sigma\!\da}

\newcommand{\Tr}{\mathrm{Tr}}
\newcommand{\Fr}{\mathrm{Fr}}

\newcommand{\tfihl}[1]{\Tr\(\Fr^s_*,\IH_{#1}\)}
\newcommand{\tfin}{\Tr\(\Fr^{-s}_*,\IH\)}

\newcommand{\scrT}{\mathcal{T}}
\newcommand{\Spec}{\mathrm{Spec}\,}
\newcommand{\Proj}{\mathrm{Proj}\,}

\newcommand{\IH}{I\! H^*}
\newcommand{\IHd}[1]{I\! H^{#1}}
\newcommand{\IHT}{I\! H^*_{\!\Td}\!}

\newcommand{\fc}{f_{{C}}}
\newcommand{\Fq}{\mathbb{F}_{\! q}}

\newcommand{\yb}{Y_{\b}}
\newcommand{\pyq}{P_Y(q)}

\newcommand{\pyqin}{P_Y(q^{-1})}
\newcommand{\pybq}{P_Y^{\b}(q)}
\newcommand{\nuybq}{\nu_{Y_{\!\b}}(q)}

\newcommand{\nuybps}{\nu_{Y_{\b}}(p^s)}

\newcommand{\paq}{P_{\!\A}(q)}
\newcommand{\DA}{\Delta_{\A}}

\newcommand{\hdq}{h_{\Delta}(q)}
\newcommand{\haq}{h_{\!\A}(q)}
\newcommand{\hsaq}{h_{\!\A}^{br}(q)}

\newcommand{\xa}{M\!(\A)}

\newcommand{\ef}{\xi_F}
\newcommand{\eb}{\xi_\b}
\newcommand{\aff}{\mathbb{A}}
\newcommand{\SR}{\mathcal{SR}}

\renewcommand{\(}{\left(}
\renewcommand{\)}{\right)}

\renewenvironment{proof}{\noindent {\bf Proof:}}{\qed \par}

\begin{document}
\baselineskip=1.1\baselineskip

\noindent {\Large \bf Intersection cohomology of hypertoric
varieties}
\bigskip\\
{\bf Nicholas Proudfoot}\footnote{Supported
by the Clay Mathematics Institute Liftoff Program
and the National Science Foundation Postdoctoral Research Fellowship.}\\
Department of Mathematics, University of Texas,
Austin, TX 78712\smallskip \\
{\bf Benjamin Webster}\footnote{Supported by the
National Science Foundation Graduate Research Fellowship.}\\
Department of Mathematics, University of California,
Berkeley, CA 94720
\bigskip
{\small
\begin{quote}
\noindent {\em Abstract.}
A hypertoric variety is a quaternionic analogue of a toric variety.
Just as the topology of toric varieties is closely related to the
combinatorics of polytopes, the topology of hypertoric varieties
interacts richly with the combinatorics of hyperplane arrangements
and matroids.  Using finite field methods, we obtain combinatorial
descriptions of the Betti numbers of hypertoric varieties, both for ordinary
cohomology in the smooth case and intersection cohomology in the singular case.
We also introduce a conjectural ring structure on the intersection cohomology
of a hypertoric variety.
\end{quote}
}
\bigskip

\noindent
Let $\Tk$ be an algebraic torus acting linearly and effectively on an affine space
$\Cn$, by which we mean a vector space over an unspecified field,
or even over the integers.  Though much of our paper is devoted to the finite
field case, for the purposes of the introduction one may simply think of a complex vector
space.
A character $\a$ of $\Tk$ defines a lift
of the action to the trivial line bundle on $\Cn$, and the corresponding
geometric invariant theory (GIT) quotient $\X = \Cn\mod_{\!\!\a}\Tk$
is a toric variety.\footnote{See \cite{P2} for an elementary
definition of GIT quotients, and an exposition of toric varieties from this perspective.
The reader with a more differential bent can think of $\X$ as a real symplectic quotient
of $\C^n$ by the compact torus sitting inside of $\Tk$.
This is not to be confused with the algebraic (or complex) symplectic quotient
of Equation \eqref{alg}.}
A hypertoric variety is a symplectic quotient
\begin{equation}\label{alg}
\M = \cot\mmod_{\!\!(\a,0)}\Tk = \minv\mod_{\!\!\a}\Tk,
\end{equation}
where $\mu:\cot\to\tkd$ is the algebraic moment map for the $\Tk$ action on $\cot$.
Over the complex numbers, this construction may be interpreted as a hyperk\"ahler
quotient \cite[\S 3]{BD}, or equivalently as a real symplectic quotient of $\minv$
by the compact form of $\Tk$.
For this reason, 
$\M$ may be thought of as a `quaternionic' or hyperk\"ahler analogue of $\X$.   
In this paper, however, we will focus on the algebro-geometric construction,
which lets us work over arbitrary fields.

The data of $\Tk$ acting on $\Cn$, along with the character $\a$,
can be conveniently encoded by an arrangement $\A$ of cooriented hyperplanes in an affine space of dimension $d=n-k$.  The topology of the
corresponding complex toric variety $\XA_\C$ 
is deeply related to the combinatorics of the
polytope cut out by $\A$ over the real numbers \cite{S1,S2}.
The hypertoric variety $\MA$ is sensitive to a different side
of the combinatorial data.  As a topological space, the complex variety $\MA_\C$ does
{\em not} depend on the coorientations of the hyperplanes \cite[2.2]{HP}, and hence has
little relationship to the polytope that controls $\XA$.
Instead, the topology of $\MA_\C$ interacts richly with the combinatorics of the matroid
associated to $\A$, as explained in \cite{Ha}.  We now describe the sort of combinatorial
structures that arise in this setting.

Let $\Delta$ be a simplicial complex of dimension $d-1$ on the ground set $\otn$.
The {\bf \boldmath$f$-vector} of $\Delta$ is the $(d+1)$-tuple $(f_0,\ldots,f_d)$,
where $f_i$ is the number of faces of $\Delta$ of cardinality $i$ (and therefore
of dimension $i-1$).  The {\bf \boldmath$h$-vector} $(h_0,\ldots,h_d)$
and {\bf \boldmath$h$-polynomial} $\hdq$ of $\Delta$
are defined by the equations $$\hdq = \sum_{i=0}^d h_i\, q^i = \sum_{i=0}^d
f_i\, q^i(1-q)^{d-i}.$$
To each simplicial complex $\D$, we associate its {\bf Stanley-Reisner ring} $\SR(\D)$,
which is defined to be the the quotient of $\C[e_1,\ldots,e_n]$ by the ideal generated
by the monomials $\prod_{i\in S}e_i$ for all non-faces $S$ of $\Delta$.
The complex $\D$ is called {\bf Cohen-Macaulay} if there exists a $d$-dimensional
subspace $L\subs\SR(\D)_1$
such that $\SR(\D)$ is a free module over the polynomial ring $\Sym L$.
Such a subspace is called a {\bf linear system of parameters}.
If $\D$ is Cohen-Macaulay and $L$ is a linear system of parameters for $\D$, then
$\SR_0(\Delta):=\SR(\Delta)\otimes_{\Sym L}\C$
has Hilbert series equal to $h_\Delta(q)$ \cite[5.9]{S3}.

Let $\A = \{H_1,\ldots,H_n\}$ be a collection of labeled hyperplanes in a vector space $V$, and let
$a_i\in V^*$ be a nonzero normal vector to $H_i$ for all $i$.
The {\bf matroid complex} $\da$ associated
to $\A$ is the collection of sets $S\subs\otn$ such that $\{a_i\mid i\in S\}$ is linearly independent.
A {\bf circuit} of $\da$ is a minimal dependent set.  Let $\sigma$ be an ordering of the
set $\otn$.  A {\bf \boldmath$\sigma$-broken circuit} of $\da$ is a set $C\smallsetminus\{i\}$,
where $C$ is a circuit, and $i$ is the $\sigma$-minimal element of $C$.
The {\bf \boldmath$\sigma$-broken circuit complex} $\bcsda$ is defined to be the collection
of subsets of $\otn$ that do not contain a $\sigma$-broken circuit.
The two complexes $\da$ and $\bcsda$ are both Cohen-Macaulay (in fact
shellable \cite[\S 7.3 $\&$ \S 7.4]{Bj});
their $h$-polynomials will be denoted $\haq$ and $\hsaq$, respectively.
As the notation suggests, the polynomial $\hsaq$ is independent
of our choice of ordering $\sigma$ \cite[\S 7.4]{Bj}.

Let $\A$ be a central hyperplane arrangement, and $\At$ a simplification of $\A$.
By this we mean that all of the hyperplanes in $\A$ pass through the origin, and $\At$
is obtained by translating those hyperplanes away from the origin in such a way so that
all nonempty intersections are generic.  Then $\MA$ is an affine cone, and
$\MAt$ is an orbifold resolution of $\MA$.
Our goal is to study the topology of the complex varieties $\MA_\C$ and $\MAt_\C$,
relating them to the combinatorics of the arrangement $\A$.
To achieve this goal, we count points on the corresponding varieties over finite fields.

Our approach to counting points on $\MAt$ is motivated by a paper of Crawley-Boevey
and Van den Bergh \cite{CBVdB}, who work in the context of representations of quivers.  
In Section \ref{smoothsec} we use an exact sequence that appeared first in \cite{CB} to
obtain a combinatorial formula for the number of $\Fq$ points of $\MAt$.
Then the Weil conjectures allow us
to translate this formula into a description of the Poincar\'e
polynomial of $\MAt_\C$ (Theorem \ref{sp}).

\begin{theorem*}
The Poincar\'e polynomial of $\MAt_\C$ coincides with the $h$-polynomial of $\DA$.
\end{theorem*}

This theorem has been proven by different means in \cite[6.7]{BD} and \cite[1.2]{HS}.  
One noteworthy
aspect of our approach is that it sheds light on a mysterious theorem of 
Buchstaber and Panov \cite[\S 8]{BP},
who produce a seemingly unrelated space with the same Poincar\'e polynomial
(see Remark \ref{buchpan}).

In the case of the singular variety $\MA$, 
we follow the example of Kazhdan and Lusztig \cite[\S 4]{KL}, who
study the singularities of Schubert varieties.  These singularities are
measured by local intersection cohomology Poincar\'e polynomials, and
Kazhdan and Lusztig obtain a recursive formula for these polynomials
using Deligne's extension of the Weil conjectures.
In our paper, we extend the argument in \cite[4.2]{KL}
to apply to more general classes of varieties.  
Roughly speaking, we consider a collection of stratified
affine cones with polynomial point count, which is closed under taking closures of strata,
and normal cones to strata.  (For details, see Theorem \ref{affine}.)
In Section \ref{strat} we give such a stratification of $\MA$, and
in Section \ref{cp} we obtain the following new result (Theorem \ref{ip}).


\begin{theorem*}
The intersection
cohomology Poincar\'e polynomial of $\MA_\C$ coincides with the $h$-polynomial
of $\bcsda$.
\end{theorem*}

Section \ref{dc} is devoted to the comparison of the two $h$-polynomials
using the decomposition theorem of \cite[6.2.5]{BBD}.  The map from $\MAt$ to $\MA$
is semismall (Corollary \ref{semismall}), hence the decomposition theorem
expresses the cohomology of $\MAt_\C$ in terms of the intersection cohomology
of the strata of $\MA_\C$ and the cohomology of the fibers of the resolution
(Equation \eqref{htdecom}).  By our previous results, we thus obtain a combinatorial
formula relating the $h$-numbers of a matroid complex to those of its
broken circuit complex.  This formula turns out to be a special case
of the Kook-Reiner-Stanton convolution formula, which is proven from a strictly
combinatorial perspective in \cite[1]{KRS}.

We note that this suggests yet another avenue leading to the computation
of the Betti numbers of $\MAt_\C$. Knowing only the intersection Betti
numbers of $\MA_\C$, we could have computed these numbers using the KRS
formula and the recursion that we obtained from the decomposition
theorem. This approach is one that generalizes naturally to other
settings.  For example, Nakajima's quiver varieties form a class of
stratified affine varieties which is closed under taking closures of
strata and normal cones to strata \cite[\S 3]{Na}.  These varieties have
semismall resolutions whose Betti numbers are relevant to the
representation theory of Ka\v c-Moody algebras, and are the subject of
an outstanding conjecture of Lusztig \cite[8]{Lu}.  If a polynomial point
count for the singular varieties could be obtained, then the
decomposition theorem would provide recursive formulas for the Betti
numbers of the smooth ones. 

Section \ref{cr} deals with the problem of ring structures.
Hausel and Sturmfels show that the
cohomology ring of $\MAt_\C$ is isomorphic to $\SR_0(\DA)$,
which strengthens Theorem \ref{sp}.
Intersection cohomology is in general only a group, so we have no analogous theorem
to prove for $\MA_\C$.  When $\A$ is a unimodular arrangement, 
however, we define a ring $R_0(\A)$ which is isomorphic
to the intersection cohomology of $\MA_\C$ as a graded vector space.
This ring does not depend on a choice of ordering $\sigma$ of the set $\{1,\ldots,n\}$,
but it degenerates flatly to $\SR_0(\bcsda)$ for any $\sigma$
(Theorem \ref{ugb}).
We conjecture that this isomorphism is natural, and that the multiplicative structure can be interpreted
in terms of the geometry of $\MA_\C$ (Conjecture \ref{nat}).

\paragraph{\bf Acknowledgments.}  The authors are grateful to 
Tom Braden, Mark Haiman, Joel Kamnitzer, Brian Osserman, 
Vic Reiner, David Speyer, and Ed Swartz for invaluable discussions.

\begin{section}{Hypertoric varieties}\label{hypertoric}
Let $\Tn$ and $\Td$ be split algebraic tori defined over $\Z$, with
Lie algebras $\tn$ and $\td$. Let $\{x_1,\ldots,x_n\}$ be a basis
for $\tn = \operatorname{Lie}(\Tn)$, 
and let $\{e_1,\ldots,e_n\}$ be the dual basis for the
dual lattice $\tnd$. Suppose given $n$ nonzero integer vectors
$\{a_1,\ldots,a_n\}\subs\td$ such that the map $\tn\to\td$ taking
$x_i$ to $a_i$ has rank $d$, and let $\tk$ be the kernel of this
map. Then we have an exact sequence
\begin{eqnarray}\label{vs}
0 \longrightarrow \tk \stackrel{\i}{\longrightarrow} \tn
\longrightarrow \td\longrightarrow 0,
\end{eqnarray}
which exponentiates to an exact sequence of groups
\begin{eqnarray}\label{tori}
0 \longrightarrow \Tk \longrightarrow \Tn
\longrightarrow \Td\longrightarrow 0,
\end{eqnarray}
where
$\Tk=\ker\!\big(\Tn\to\Td\big)$. Thus $\Tk$ is an algebraic group with
Lie algebra $\tk$, which is connected if and only if the vectors
$\{a_i\}$ span the lattice $\td$ over the integers. Every algebraic
subgroup of $\Tn$ arises in this way.

Consider the cotangent bundle $\cot\cong\Cn\times\Cnd$ along with its natural
algebraic symplectic form $$\omega = \sum dz_i \wedge dw_i,$$
where $z$ and $w$ are coordinates on $\Cn$ and $\Cnd$, respectively.
The restriction to $\Tk$ of the standard action of
$\Tn$ on $\cot$ is hamiltonian, with moment map
$$\muc(z,w) = \i^*\sum_{i=1}^n (z_i w_i)\, e_i.$$
Suppose given an integral element $\a\in\tkd$.  This descends via the 
exponential map to a 
character of $T^k$, which defines a 
lift of the action of $\Tk$
to the trivial bundle on $\cot$.
The symplectic quotient
$$\Ma=\cot\mmod_{\!\!\a}\Tk =
\mu^{-1}(0)\mod_{\!\!\a}\Tk$$
is called a {\bf hypertoric variety}.  
Here the second quotient is a projective GIT quotient\footnote{For a careful
treatment of geometric invariant theory over the integers, see Appendix B
of \cite{CBVdB}.}
$$\mu^{-1}(0)\mod_{\!\!\a}\Tk := \Proj\bigoplus_{m=0}^{\infty}
\Big\{f\in\mathcal{O}_{\!\minv}\bigmid \nu^*(f) = \a^m\otimes f\hs\Big\},$$
where $$\nu^*:\O_{\mu^{-1}(0)}\to\O_{T^k\times\mu^{-1}(0)}\cong\O_{T^k}\otimes\O_{\mu^{-1}(0)}$$
is the map on functions induced by the action map $\nu:T^k\times\mu^{-1}(0)\to\mu^{-1}(0)$.
If $\a$ is omitted from the subscript, it will be understood to be equal to zero.
The hypertoric variety $\Ma$ is a symplectic variety of dimension $2d$, and
admits an effective hyperhamiltonian action of the torus $\Td=\Tn/\Tk$,
with moment map
$$\Phi[z,w] = \sum_{i=1}^n (z_i w_i)\,e_i
\in\ker(\i^*) = \tdd.$$
Here $[z,w]$ is used to denote the image in $\Ma$ of a pair 
$(z,w)$ with closed $\Tk$-orbit in $\mu^{-1}(0)$.

\begin{remark}
The word `hypertoric' comes from the fact that the complex variety $\Ma_\C$
may be constructed as a hyperk\"ahler quotient of $T^*\C^n$ by the compact
real form of $\Tk$, thus making it a `hyperk\"ahler analogue'
of the toric variety $\X = \C^n\mod_{\!\!\a}\Tk$.
This was the original approach of Bielawski and Dancer {\em\cite{BD}}, who
used the name `toric hyperk\"ahler manifolds'.  For more on the general
theory of hyperk\"ahler analogues of K\"ahler quotients, see {\em\cite{P1}}.
\end{remark}

It is convenient to encode the data that were used to construct
$\Ma$ in terms of an arrangement of affine hyperplanes in
$\tdd$,
with some additional structure. A {\bf weighted,
cooriented, affine hyperplane} $H\subs\tdd$ is a
hyperplane along with a choice of nonzero integer normal vector
$a\in\td$.  Here ``affine'' means that $H$ need not pass through the origin, and
``weighted'' means that $a$ is not required to be primitive. Let
$r = (r_1,\ldots,r_n)\in\tnd$ be a lift of $\a$ along $\i^*$, and
let $$H_i =\{v\in\tdd \mid v\cdot a_i + r_i = 0\}$$ be the
weighted, cooriented, affine hyperplane with normal vector
$a_i\in\td$. We will denote the arrangement $\{H_1,\ldots,H_n\}$
by $\A$, and the associated hypertoric variety by $\MA$.
Choosing a different lift $r'$ of $\a$ corresponds geometrically to translating
$\A$ inside of $\tdd$.
The rank of the lattice spanned by the vectors $\{a_i\}$ is called the {\bf rank} of $\A$;
in our case, we have already made the assumption that $\A$ has rank $d$.
Observe that $\A$ is defined over the integers, and therefore may be realized
over any field.  Intuitively, it is useful to think of $\A$ as an arrangement
of real hyperplanes, but in Section \ref{cp} we will need to consider the
complement of $\A$ over a finite field.

\begin{remark}
We note that we allow repetitions of hyperplanes in our
arrangement ($\A$ may be a multi-set), and that a repeated
occurrence of a particular hyperplane is {\em not} the same as a
single occurrence of that hyperplane with weight 2. On the other
hand, little is lost by restricting one's attention to
arrangements of distinct hyperplanes of weight one.
\end{remark}

We call the arrangement $\A$ {\bf simple} if every subset of $m$
hyperplanes with nonempty intersection intersects in codimension
$m$. We call $\A$ {\bf unimodular} if
every collection of $d$ linearly independent vectors
$\{a_{i_1},\ldots,a_{i_d}\}$ spans $\td$ over the integers.
An arrangement which is both simple and unimodular is called {\bf smooth}.

\begin{theorem}\label{simple}{\em\cite[3.2 $\&$ 3.3]{BD}}
The hypertoric variety $\MA$ has at worst orbifold (finite quotient)
singularities if and only if $\A$ is simple,
and is smooth if and only if $\A$ is smooth.
\end{theorem}

Let $\A=\{H_1,\ldots,H_n\}$ be a
{\bf central} arrangement, meaning that $r_i = 0$ for all $i$.
Let $\At = \{\tilde
H_1,\ldots, \tilde H_n\}$ be a {\bf simplification} of $\A$, by
which we mean an arrangement defined by the same vectors
$\{a_i\}\subset\td$, but with a different choice of $r\in\tnd$, such that
$\At$ is simple. This corresponds to translating each of the
hyperplanes in $\A$ away from the origin by some generic amount.
We then have
\begin{eqnarray*}
\MA \,=\, \Proj\bigoplus_{m=0}^\infty\mathcal{O}_{\!\minv}^{\Tk}
\,=\, \Spec \mathcal{O}_{\!\minv}^{\Tk}
\,=\, \Spec \mathcal{O}_{\MAt},
\end{eqnarray*}
hence there is a surjective, projective map $\pi:\MAt\to\MA$.
Geometrically, $\pi$ may be understood to be the map induced
by the $\Tk$-equivariant
inclusion of $\mcinvast$ into $\mcinv$, where $\mcinvast$ is the
stable locus for the linearization of the $\Tk$ action given by
$\a = \i^*(r)$.
The central fiber $\LAt = \pinv(0)$ is called the {\bf core} of
$\MAt$.

\begin{theorem}\label{core}{\em\cite[\S 6]{BD},\cite[6.4]{HS}}
The core $\LAt$ is isomorphic to a union of toric varieties with
moment polytopes given by the bounded complex of $\At$.
Over the complex numbers,
$\LAt_\C$ is a $\Td_\R$-equivariant deformation retract of
$\MAt_\C$, where $\Td_\R\cong U(1)^d$ is the compact
real form of $\Tdc$.
\end{theorem}

It follows that the dimension of the core is at most $d$, with equality
if and only $\A$ is coloop-free (see Remark \ref{sl}).

\begin{example}\label{ex}
The two arrangements pictured below are each simplifications
of a central arrangement of four hyperplanes in $\R^2$,
in which the second and third hyperplane coincide.  All hyperplanes are taken with
weight 1, and coorientations may be chosen arbitrarily.

\begin{figure}[h]
  \centerline{\epsfig{figure=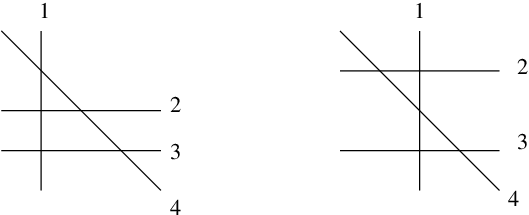, height=2cm}}
\end{figure}
\noindent Consider the complex hypertoric varieties associated to these two
arrangements.  Both are
obtained as symplectic quotients of $T^*\mathbb C^4$ by the same $T^2$
action, but with different choices of character.  Both varieties
are resolutions of the affine variety given by the associated
central arrangement. The hypertoric variety associated to the
left-hand arrangement has a core consisting of a projective plane
glued to a Hirzebruch surface along a projective line.  The hypertoric
variety associated to the right-hand arrangement has a core
consisting of two projective planes glued together at a point.
As manifolds, they are
diffeomorphic, as are any two complex hypertoric varieties
corresponding to different simplifications of the same central
arrangement {\em\cite[2.1]{HP}}.
\end{example}
\end{section}

\begin{section}{The stratification}\label{strat}
Let $\A$ be a rank $d$ central arrangement of $n$ weighted, cooriented hyperplanes in $\tdd$.
Our goal for this section is to define and analyze a stratification
of the singular affine 
variety $\MA$.  This stratification will be a refinement
of the Sjamaar-Lerman stratification, introduced for real symplectic quotients
in \cite{SL}, and adapted to the algebraic setting in \cite[\S 3]{Na}.
Our refinement will prove to be more natural from a combinatorial perspective
(see Remark \ref{sl}).

Given any subset $S\subs\otn$, let $H_S = \cap_{i\in S}H_i$. A
{\bf flat} of $\A$ is a subset $F\subs\otn$ such that $F = \{i\mid
H_i\sups H_F\}$.  We let $L(\A)$ denote the lattice of flats for
the arrangement $\A$.  For any flat $F$, we define the {\bf
restriction}
$$\AF := \{H_i\cap H_F\mid i\notin F\},$$
an arrangement of $|F^c|$ hyperplanes in the affine space $H_F$,
and the {\bf localization}
$$\Af := \{H_i/H_F\mid i\in F\},$$
an arrangement of $|F|$ hyperplanes in the affine space
$\tdd/H_F$. The lattice $L(\AF)$ is isomorphic to the sublattice
of $L(\A)$ consisting of those flats which contain $F$; likewise,
$L(\Af)$ may be identified with the sublattice of $L(\A)$ consisting of flats
contained in $F$. We define the {\bf rank} of a flat $\rk F =
\rk\Af$, and the {\bf corank} $\crk F = \rk\AF = \rk\A - \rk F$.
Given a simplification $\tilde\A$ of $\A$, there is a natural
simplification $\tilde\Af$ of the localization.

We now fix notation regarding the various tori associated to the
localization and restriction of $\A$ at $F$.
The $\Af$ analogue of the exact sequence \eqref{vs} is
$$0\to\hat{\mathfrak{t}}\to\tf\to\mathfrak{t}^{\rk F}\to 0,$$
where $\tf$ is the coordinate subtorus of $\tn$ supported on $F$,
$\mathfrak{t}^{\rk F}\cong H_F$ is the image of $\tf$ in $\td$,
and $\hat{\mathfrak{t}} = \tk\cap\tf$.  Similarly, the restriction $\AF$
corresponds to an exact sequence
$$0\to\mathfrak{t}\to\tfc\to\mathfrak{t}^{\crk F}\to 0,$$
where $\tfc$ is the coordinate subtorus of $\tn$ supported on $F^c$,
$\mathfrak{t}^{\crk F}\cong \td/H_F$,
and $\mathfrak{t} = \tk/\hat{\mathfrak{t}}$.
The tori $\hat{T}, T^{F}, T^{\rk F}$ and $T,T^{F^c},\Tcrkf$
are defined analogously, as in the exact sequence \eqref{tori}.
Let $T^*\!\mathbb A^{\! F}$ and $T^*\!\mathbb A^{\! F^c}$ be the
cotangent bundles of the corresponding coordinate subspaces of $A^n$.
Then the hypertoric varieties $\MAf$ and $\MAF$ are obtained as symplectic
quotients of $T^*\!\mathbb A^{\! F}$ and $T^*\!\mathbb A^{\! F^c}$
by $\hat T$ and $T$, respectively.

\begin{proposition}\label{ressub}
Let $F$ be a flat of $\A$.  The subvariety of $\MA$
given by the equations $z_i=w_i=0$ for all $i\in F$
is isomorphic to $\MAF$.
\end{proposition}

\begin{proof}
The inclusion of
$T^*\!\mathbb A^{\! F^c}$ into $\cot$ is $\Tk$-equivariant, where the action
of $\Tk$ on $T^*\!\mathbb A^{\! F^c}$ factors through $T$. The Lie coalgebra
$\mathfrak{t}^*$ of $T$ includes into $\tkd$, and the $\Tk$-moment
map on $\cot$ restricts to the $T$-moment map on $T^*\!\mathbb A^{\! F^c}$, as shown in the
following diagram.
$$\xymatrix{\cot \ar[r]^(.4){\quad\mu}& \tkd\\
T^*\!\mathbb A^{\! F^c}\ar[r]^{\quad\mu^F}\ar[u] &\ar[u]
\mathfrak{t}^*}$$
Thus we have
$$\MAF = T^*\!\mathbb A^{\! F^c} \mmod T = \left({\mu^{F}}\right)^{-1}\!(0)\bigmod T
\cong \Big(T^*\!\mathbb A^{\! F^c} \cap \minv\Big)\bigmod\Tk,$$
which is cut out of $\MA=\minv/\Tk$
by the equations $z_i=w_i=0$ for all $i\in F$.
\end{proof}

\begin{remark}
The subvariety $\MAF$ in Proposition \ref{ressub} may also be described as the preimage of $H_F^\R\oplus H_F^\C$
along the hyperk\"ahler moment map; see the proof of Proposition \ref{cohomdecom}.\footnote{The statement of Remark 2.2
appearing in the published version of this paper is incorrect; we thank Linus Setiabrata for pointing this out.}
\end{remark}

Let $S,Y$, and $Z$ be schemes over an arbitary ring, with $Z\subs Y$ a locally closed subscheme $Z$, and $s\in S$ a basepoint.  
We will say that $S$ is a {\bf normal slice}
to $Z$ in $Y$ if there exists a collection of \'etale open cover $\{U_\alpha\}$ 
of a neighborhood of $Z$, such that each $U_\alpha$
admits a dominant \'etale map to $S\times\mathbb A^{\dim Z}$,
with $U_\alpha\cap Z$ dominating $\{s\}\times\mathbb A^{\dim Z}$.
Since \'etale maps are locally invertible in the analytic category, this
implies that our definition of normal slices for complex schemes agree with the usual
notion.
In other words,
there is an analytic neighborhood of $Z$ in $Y$ which is
locally biholomorphic to a neighborhood of $Z\times\{s\}$ 
in $Z\times S$.

A {\bf stratification} of a scheme $Y$ over $\Z$
is a partition of $Y$ into smooth, locally closed subschemes $Y_\b$
indexed by a finite poset $B$, 
along with normal slices $S_\b$,
with the property that for all $\b\in B$,
$$\overline{Y_\b}=\bigsqcup_{\b'\leq\b}Y_{\b'}.$$
To define a stratification of $\MA$, we begin by putting
$$\MAopen = \Big\{[z,w]\in\MA\bigmid 
\text{there does not exist $i$ such that } z_i=w_i=0\Big\}.$$
The identification of $\MAF$ with a subvariety of $\MA$ as in Proposition
\ref{ressub} induces the identification
\begin{equation}\label{stratum}
\MAFopen = \Big\{[z,w]\in\MA\bigmid
z_i=w_i=0\text{ if and only if }i\in F\Big\}.
\end{equation}
By \cite[3.1]{BD}, any point $[z,w]\in\MA$ has the property that
the set $\{i\mid z_i=w_i=0\}$ is a flat, hence we have a decomposition
\begin{equation}\label{decomp}
\MA = \displaystyle\bigsqcup_F\MAFopen.
\end{equation}
One interpretation of our decomposition is that points are grouped
according to the
stabilizers in $\Tn$ of their lifts to $\minv$.  This is therefore a 
refinement of the Sjamaar-Lerman stratification
\cite[\S 2]{SL}, which groups points by their stabilizers in
the subtorus $\Tk$.  
It follows from the defining property of a moment map that $\MAopen$
is smooth, and therefore that each piece of the decomposition is smooth.
To see that our
decomposition is a stratification we must produce normal slices to the strata,
which we will do in Lemma \ref{etneigh}.
The largest stratum $\MAopen$ will be referred to as
the {\bf generic stratum} of $\MA$.

\begin{remark}\label{sl}
If $[z,w]\in\MAFopen$, then the stabilizer of $(z,w)\in\minv$ is equal
to $\hat{T} = \Tk\cap T^F$.
An element $i\in F$ is called a {\bf coloop} of $F$ if $a_i$ may not be expressed as
a linear combination of $\{a_j\mid j\in F\smallsetminus\{i\}\}$.
If $F$ and $G$ are two flats, then $\Tk\cap\TF = \Tk\cap\TG$ if and only
if $F$ and $G$ agree after deleting all coloops.  
Hence the Sjamaar-Lerman stratification
of $\MA$ is naturally indexed by coloop-free flats, rather than all flats.
\end{remark}


\begin{lemma}\label{etneigh}
The variety $\MAf$ is a normal slice to  $\MAFopen\subs\MA$,
thus the decomposition of Equation \eqref{decomp} is a stratification.
\end{lemma}

\begin{proof}
Let $S$ be a subset of $F^c$ such that the coordinate vectors
in $(\tfc)^*$ indexed by $S$ descend to a basis of $\mathfrak{t}^*$, and
let 
\begin{equation}\label{ws}
W_S = \Big\{[z,w]\in\MA\bigmid z_i\neq 0\text{ for all $i\in S$}\Big\}.
\end{equation}
Any element of $W_S$ may be represented by an element $(z,w)\in\minv$
such that $z_i = 1$ for all $i\in S$, and any two such representations
differ by a element of the subtorus $\hat{T}\subs\Tk$.
Next observe that the coordinate projection of $\cot$ onto $T^*\!\mathbb A^{\! F}$
takes $\minv$ to the zero set of the moment map for the action of $\hat{T}$
on $T^*\!\mathbb A^{\! F}$.
We may therefore define a map $p_F:W_S\to\MAf$ by taking an element of $W_S$,
representing it in the form described above, and projecting to the $F$ 
coordinates. This map is smooth on the locus $\MAFopen\cap W_S$.

Suppose that $y\in W_S$.
By \cite[2.2.14]{BLR}, there is a neighborhood $U$
of $y$ in $W_S$ and a smooth map $\eta:U\to\aff^{\crk F}$
such that the restriction of $p_F$ to $\eta^{-1}(0)$ is \'etale.
Let $$\vartheta=\eta\times p_F:U\to\aff^{\crk F}\times\MAf.$$
Then the derivative of $\vartheta$ at $y$ is a surjection, 
hence $\vartheta$ is smooth at $y$.
Since its source and target have the same dimension, it must be \'etale.

If $y$ is not in $S$, then we may modify the definition of $W_S$ by changing
some of the $z_i$ to $w_i$ in Equation \eqref{ws}, and adjust the definition of the map
$p_F$ accordingly.  Then $y$ will be contained in the new set $W_S$, and the proof
will go through as before.
\end{proof}

\vspace{\baselineskip}
We next prove a result similar to Lemma \ref{etneigh} by working purely in the
analytic category.  The advantage of Lemma \ref{slice} is that we obtain a 
statement that is compatible with the affinization map $\pi$, which will be useful
in Section \ref{dc}.

\begin{lemma}\label{slice}
For all $y\in\MAFopen$ there is an analytic neighborhood $\mathfrak{U}$ of $y\in\MA_\C$
and a map $\varphi:\mathfrak{U}\to \MAFopen_\C\times\MAf_\C$
such that $\varphi(y) = (y,0)$, and $\varphi$ is a diffeomorphism
onto its image.  Furthermore, there is a map
$\tilde\varphi:\pi^{-1}\!\(\mathfrak{U}\)\to\MAFopen_\C\times\MAft_\C$
which covers $\varphi$, and is also a diffeomorphism onto its image.

$$\xymatrix{\MAt_\C \ar[d]^{\pi}& \pi^{-1}(\mathfrak{U})\ar[d]^\pi\ar[l]\ar[r]^(.3){\tilde\varphi}
& \MAFopen_\C \times\MAft_\C\ar[d]^{\mathrm{id}\times\pi_F}\\
\MA_\C & \mathfrak{U} \ar[l]\ar[r]^(.3){\varphi}&
\MAFopen_\C\times\MAf_\C}$$
\end{lemma}

\begin{proof}
Let $\bar y\in T^*\C^n$ be a representative of $y$.
Since $y$ is contained in the stratum $\MAFopen_{\C}$, we may assume that the $\ith$
coordinates of $\bar y$ are zero for all $i\in F$.
Let $V = T_{\bar y}\big(\Tkc\cdot\bar y\big)$ 
be the tangent space to the orbit of $\Tk_\C$ through $\bar y$.
Then $V \subs T^*\C^n$ is isotropic with respect to the symplectic form $\omega$,
and the inclusion of $T^*\mathbb C^{F}$ into $T^*\C^n$ induces a $\Thc$-equivariant
inclusion of $T^*\mathbb C^{F}$ into the quotient $V^{\omega}\!/V$, 
where $V^{\omega}$ is the symplectic perpendicular space to $V$
inside of $T^*\C^n$.  The torus $\Thc$ acts trivially 
on the quotient of $V^{\omega}\!/V$ by 
$T^*\mathbb C^{F}$, which may be identified with the tangent space $T_y\MAFopen_\C$.
The lemma then follows from the discussion in \cite[\S 3.2 $\&$ 3.3]{Na}.
\end{proof}

\begin{corollary}\label{loctriv}
The restriction of $\pi$ to $\pi^{-1}\Big(\MAFopen_\C\Big)$ is a locally
trivial topological fiber bundle over the stratum $\MAFopen_\C$,
with fiber isomorphic to the core $\mathcal{L}(\Aft)_\C\subs\MAft_\C$.
\end{corollary}

If $Y = \sqcup_{\b\in B}\yb$ is a stratified space and $f:X\to Y$ is a map,
then $f$ is called {\bf semismall} if for all $y_\b\in\yb$, the dimension
of $f^{-1}(y_\b)$ is at most half of the codimension of $\yb$ in $Y$.  This seemingly
arbitrary condition can be motivated by the observation that
$$f\hs\text{  semismall}\hs\hs\iff\hs\hs\dim X=\dim Y = \dim \Big(Y\times_X Y\Big).$$

\begin{corollary}\label{semismall}
The map $\pi$ is semismall.
\end{corollary}

\begin{proof}
$2\dim\mathcal{L}(\Aft) \leq \dim\MAft = \codim\MAFopen$, with equality if
and only if $F$ is coloop-free.
\end{proof}

\begin{remark}\label{ok}
The decomposition defined by Equations \eqref{stratum} and \eqref{decomp}
make sense for noncentral arrangements as well as central ones.
For most of the paper we will
use this decomposition only in the central case, 
but in Section \ref{cr} we will consider the generic stratum of an arbitrary
hypertoric variety.
\end{remark}
\end{section}

\begin{section}{The Betti numbers of \boldmath$\MAt$}\label{smoothsec}
Let $X$ be a variety defined over the integers, and let $q$
be a prime power.
By an {\bf \boldmath$\Fq$ point} of $X$, we mean a closed point of the
variety $X_{\Fq} = X\otimes_{\Z}\Fq$.
We say that
$X$ has {\bf polynomial point count} if there exists a polynomial
$\nu_X(q)$ such that, when $q$ is a power of a sufficiently large
prime, $\nu_X(q)$ is equal to the number of $\Fq$ points of $X$.
For the remainder of the paper, when we refer to the number of
$\Fq$ points on a given variety, we will always implicitly
assume that $q$ is a power of a sufficiently large prime.
Suppose that $X$ has at worst orbifold singularities, so that the cohomology
of $X_\C$ is Poincar\'e dual to the compactly supported cohomology.
The Betti numbers for the compactly supported cohomology of $X_\C$ agree with the Betti
numbers for compactly supported $\ell$-adic \'etale cohomology of 
$X_{\Fpbar}$ for large enough primes $p$ \cite[6.1.9]{BBD}.
If the $\ell$-adic \'etale cohomology of $X_{\Fpbar}$ 
is pure, we may use the Lefschetz fixed point theorem in $\ell$-adic \'etale cohomology
to relate the Betti numbers of $X_\C$ to the number of $\Fq$-points on $X$ 
(see for example \cite[A.1]{CBVdB}).

\begin{theorem}\label{smoothweil}
Suppose that $X$ has polynomial point count
and at worst orbifold singularities, 
and that the $\ell$-adic \'etale cohomology of the variety $X_{\Fqbar}$ 
is pure. 
Then $X_{\C}$ has Poincar\'e polynomial 
$P_X(q) = q^{\dim X}\cdot\nu_X(q^{-1})$, where $q$ has degree 2.
(In particular, the odd cohomology of $X$ vanishes.)
\end{theorem}


The purpose of this section is to apply Theorem \ref{smoothweil}
to compute the Poincar\'e polynomial $P_{\tilde A}(q)$ of $\MAt_\C$.
The fact that the $\ell$-adic \'etale cohomology of $\MAt_{\Fqbar}$ is pure
follows from \cite[A.2]{CBVdB},\footnote{The authors state this theorem only
for smooth varieties, but their argument clearly extends to the orbifold case.}
using the $\mathbb{G}_m$-action
studied in \cite{HP}.

Recall that $\MAt$ is defined as the GIT quotient
$\mu^{-1}(0)\mod_{\!\!\a}\Tk$, where $\mu$
is the moment
map for the action of $\Tk$ on $\cot$. For any $\lambda\in\tkd$, let
$$\MAtl = \mu^{-1}(\lambda)\mod_{\!\!\a}\Tk.$$ If $\lambda$ is a regular
value of $\mu$, then $\Tk$ will act locally freely on $\mil$, 
meaning that the stabilizer in $\Tk$ of any point in $\mil$
is finite over any field.\footnote{When we speak of a finite
subgroup of a torus that is defined over a finite field,
we always mean that the subgroup remains finite after passing to the
algebraic closure.}
This in turn implies that the
GIT quotient of $\mil$ by $\Tk$ over any algebraically closed field
will be an honest geometric quotient. 
Fix a regular value $\lambda$.
By an argument completely analogous to that of Nakajima's appendix to \cite{CBVdB},
the varieties $\MAt$ and $\MAtl$ have the same point count.
Thus Theorem \ref{smoothweil} tells us
that we can compute $P_{\tilde A}(q)$ by counting points on $\MAtl$
over finite fields.

\begin{lemma}\label{locfree}
Let $X$ be a variety defined over $\Fq$, let $T$ be a split torus of
rank $k$ acting on $X$, and let $T'$ be a (possibly disconnected)
rank $\ell$ subgroup 
of $T$ which acts locally freely.  Then the number of $\Fq$ points
of $X$ is equal to $(q-1)^\ell$ times the number of $\Fq$ points of $X/T'$.
\end{lemma}

\begin{proof}
By Hilbert's Theorem 90, every $T$-orbit in $X$ which defines an $\Fq$ point
of $X/T$ contains an $\Fq$ point of $X$.  This tells us that we may count 
points on $X$ $T$-orbit by $T$-orbit, and thereby reduce to the case where 
$T$ acts transitively. In this case $X$ is isomorphic to a split torus, and 
$T'$ acts on $X$ via a homomorphism with finite kernel.  
Thus $X/T'$ is a split torus with $\rk X/T'=\rk X-\ell$.  This completes the proof.
\end{proof}





\begin{corollary}\label{milcount}
The number of $\Fq$ points of $\mil$ is equal
to $(q-1)^k$ times the number of $\Fq$ points of $\MAtl$.
\end{corollary}


\begin{proposition}\label{smoothcount}
The variety $\MAtl$ has polynomial point count, with
$$\nu_{\MAtl}(q) = q^{2d}\cdot h_{\A}(q^{-1}).$$
\end{proposition}

\begin{proof}
For any element $z\in\Cn$, we have 
an exact sequence\footnote{The analogous exact sequence 
in the context of representations of quivers
first appeared in \cite[3.3]{CB}, and was used to count points on quiver varieties
over finite fields in \cite[\S 2.2]{CBVdB}.}
\begin{equation}\label{cb}
0\to\{w\mid\mu(z,w)=0\}\to T_z^*\Cn\overset{\mu(z,-)}{\longrightarrow}
\tkd\to\stab(z)^*\to 0,
\end{equation}
where $\stab(z)^* = \tkd/\stab(z)^\perp$ is the Lie coalgebra of the
stabilizer of $z$ in $\Tk$.
Consider the map $\phi:\mil\to\Cn$
given by projection onto the first coordinate.
By exactness of \eqref{cb} at $\tkd$, we have
$$\Im\phi = \{z\mid\lambda\cdot\stab(z) = 0\}
= \{z\mid\stab(z) = 0\},$$
where the last equality follows from the fact that $\lambda$ is a regular value.
Furthermore, we see that for $z\in\Im\phi$, $\phi^{-1}(z)$
is a torsor for the $d$-dimensional vector space
$\{w\mid\mu(z,w)=0\}$.  Hence the number of $\Fq$ points of $\mil$
is equal to $q^d$ times the number of $\Fq$ points of $\Cn$ at which
$\Tk$ is acting locally freely.

A point $z\in\Cn$ is acted upon locally freely by $\Tk$ if and
only if $\{i\mid z_i = 0\}\in\DA$.  Hence the total number of such points
over $\Fq$ is equal to
\begin{eqnarray*}
\sum_{S\in\DA}(q-1)^{n-|S|} &=& \sum_{i=0}^d f_i(\DA)\cdot(q-1)^{n-i}\\
&=& (q-1)^k\sum_{i=0}^d f_i(\DA)\cdot(q-1)^{d-i}\\
&=& (q-1)^k\cdot q^d\cdot\sum_{i=0}^d f_i(\DA)\cdot q^{-i}(1-q^{-1})^{d-i}\\
&=& (q-1)^k\cdot q^d\cdot h_{\A}(q^{-1}).
\end{eqnarray*}
To find the number of $\Fq$ points of $\MAtl$
we multiply by $q^d$ and divide by $(q-1)^k$,
and thus obtain the desired result.
\end{proof}

\vspace{\baselineskip}
Theorem \ref{smoothweil} and Proposition \ref{smoothcount}, along with the observation
that $\nu_{\MAt}(q) = \nu_{\MAtl}(q)$, combine to give
us the Poincar\'e polynomial of $\MAt$.

\begin{theorem}\label{sp}
The Poincar\'e polynomial of $\MAt_\C$ coincides with the $h$-polynomial
of the matroid complex associated to $\A$, that is 
$P_{\tilde\A}(q) = \haq$.
\end{theorem}

\begin{remark}\label{buchpan}
Implicit in the work of Buchstaber and Panov {\em\cite[\S 8]{BP}} is a calculation of the
cohomology ring of the nonseparated complex variety
$W/\Tk_\C$, where $W\subs\C^n$ is the locus of points at which $\Tk_\C$ acts locally freely.
Their description of this ring coincides with the description of $H^*(\MAt_\C)$
that we will give in Theorem \ref{smooth}, due originally to {\em\cite{Ko,HS}}.  
We now have an explanation of why these rings are the same:  $\MAt_\C$
is homeomorphic to $\MAtl_\C$, which, by the proof of Proposition
\ref{smoothcount}, is an affine space bundle
over $W/\Tk_\C$.
\end{remark}

\end{section}

\begin{section}{The Betti numbers of \boldmath$\MA$}\label{cp}
Our aim in this section is to prove an analogue of Theorem \ref{sp}
for the intersection cohomology of the singular variety $\MA$.
Intersection cohomology was defined for complex varieties in
\cite{GM1,GM2}.  The sheaf theoretic definition naturally extends to
an $\ell$-adic \'etale version for varieties in positive characteristic,
which was studied extensively in \cite{BBD}.
  
Let $Y$ be a variety of dimension $m$, defined over the integers, with
a stratification
$$Y = \displaystyle\bigsqcup_{\b\in B}\yb.$$
Let us suppose further that for all $\b$, the normal slice $S_\b$ to the
stratum $Y_\b$ is an affine cone,
meaning that it is equipped with an action of the multiplicative
group $\mathbb G_m$ having the basepoint $s$ as its unique fixed point,
and that $s$ is an attracting fixed point.
Let $\IH(Y)$ denote the global $\ell$-adic \'etale intersection
cohomology of $Y_{\Fpbar}$ for $p$ a large prime, and $\IH_{\!\b}(Y)$ the local 
$\ell$-adic \'etale intersection cohomology
at any point in the stratum $(Y_\b)_{\Fpbar}$.
Since local intersection cohomology is preserved by any \'etale map, 
and the global intersection cohomology of a cone is the same as the local
intersection cohomology at the vertex by \cite[\S 3]{KL}, 
we have natural isomorphisms
\begin{equation}\label{glo-loc-iso}
\IH_{\!\b}(Y)\cong \IH_{\!s}(S_\b)\cong\IH(S_\b)
\end{equation}
for all $\b\in B$.

In this case,
let $$\pyq = \sum_{i=0}^{m-1}\dim I\! H^{2i}(Y)\cdot q^i$$
be the even degree intersection cohomology Poincar\'e polynomial of $Y\!$,
and let $$\pybq = \sum_i\dim I\! H^{2i}_{\b}(Y)\cdot q^i$$
be the corresponding Poincar\'e polynomial for the local
intersection cohomology at a point in $\yb$.
(In the cases of interest to us, odd degree global and local cohomology
will always vanish.)
Provided that $p$ is chosen large enough,
these polynomials agree with the Poincar\'e polynomials
for global and local topological intersection cohomology
of the complex analytic space $Y_\C$ by \cite[6.1.9]{BBD}.

Let $\scrT$ be a class of stratified schemes over $\Z$ 
satisfying the following
two conditions.
\begin{enumerate}\renewcommand{\labelenumi}{(\arabic{enumi})}
 \item For each stratum $Y_\b$ of $Y\in\scrT$, the normal slice $S_\b$ to $Y_\b$ in $Y$
 is isomorphic to an element of $\scrT$.
\item For each $Y\in\scrT$, the 
group $\IH(Y)$ is pure.
\end{enumerate}
The following analogue of Theorem \ref{smoothweil} is a generalization 
of the main result of \cite[\S 4]{KL}, in which the place of $\scrT$ was taken 
by the class consisting of the intersections of Schubert varieties and
opposite Schubert cells.

\begin{theorem}\label{affine} 
Suppose that 
every element of $\scrT$ has polynomial point count. 
Then all global and local intersection cohomology groups of elements of $\scrT$ 
vanish in odd degree,
and 
for all $Y\in \scrT$, we have
\begin{equation}\label{weilrecur}
q^m\cdot\pyqin = \displaystyle\sum_{\b\in B}\pybq\cdot\nuybq.
\end{equation}
\end{theorem}

\begin{proof}
  Let $\Fr_*^s:\IH(Y_{\Fpbar})\to\IH(Y_{\Fpbar})$ be the map induced by
  the $s^\text{th}$ power of the Frobenius automorphism
  $\Fr:Y_{\Fpbar}\to Y_{\Fpbar}$.  Purity of $\IH$ implies that the
  eigenvalues of $\Fr^s_*$ on $\IHd {i}$ all have absolute value
  $p^{is/2}$.  The polynomial point count hypothesis implies that each
  eigenvalue $\alpha$ of $\Fr_*$ must satisfy $\alpha^s=f(p^s)$ for some
  polynomial $f$.  This is only possible when $f(x)=x^{i/2}$ and $i$ is
  even.  Thus odd cohomology vanishes, and the eigenvalues of $\Fr_*^s$
  on $\IHd {2i}$ is $p^{si}$.  Since $\IH_\b$ is isomorphic to the
  global intersection cohomology of the normal slice $S_\b$, and
  $\IH(S_\b)$ is pure by conditions (1) and (2), the odd cohomology
  vanishes and eigenvalues of $\Fr_*^s$ are all $p^{is}$ for $\IH_\b$ as
  well.  Thus
\begin{equation}\label{frobei}
P_{Y}(p^s)=\Tr\(\Fr_*^s,\IH\) \quad\text{and}\quad P^\b_Y(p^s)=\Tr\(\Fr_*^s,\IH_\b\).
\end{equation}

By Poincar\'e duality and the Lefschetz formula
\cite[II.7.3 \& III.12.1(4)]
{KW}, we have
\begin{align}\label{frobrecur}
p^{ms}\cdot\tfin&=\tfihl{c} \nonumber\\ 
 &=\sum_{\Fr^s(y)=y}\tfihl{y} \\ 
 &=\sum_{\beta\in B} \nuybps \cdot \tfihl{\!\b} \nonumber.
\end{align}

Equation \eqref{weilrecur} follows immediately from substitution of
\eqref{frobei} into \eqref{frobrecur}.
\end{proof}


\vspace{\baselineskip}
Let $\A$ be a central hyperplane arrangement as in Section \ref{strat},
and let $\paq = P_{\MA}(q)$.
Our goal is to show that $\MA$ has polynomial point count, and
to use Theorem \ref{affine} to compute its 
intersection cohomology Poincar\'e polynomial.
Let $$\xa = \tdd\smallsetminus\bigcup_{i=1}^n H_i$$ be the complement
of $\A$ in $\tdd$.  Then $\xa$ has polynomial point count,
and the polynomial $\chi_{\A}(q) := \nu_{\xa}(q)$ is
known as the {\bf characteristic polynomial} of $\A$ \cite[2.2]{At}.
Let $r(\A)$ denote the number
of components of the real manifold $\xa_\R$.

\begin{proposition}\label{count}
The hypertoric variety $\MAopen$ has polynomial point count, with
$$\nu_{\becircled\MA}(q)=
(q-1)^d\cdot\displaystyle\sum_F\chi_{\AF}(q)\cdot r(\Af).$$
\end{proposition}

\begin{proof}
Consider the decomposition
$$\tdd = \bigsqcup_F M\!(\AF)$$
into complements of restrictions of the arrangement $\A$ to various flats.
We will count points of $\MAopen$ on the individual fibers of 
the moment map $\Phi:\MA\to\tdd$,
and then add up the contributions of each fiber.
The fiber $$\Phi^{-1}(0) = \big\{[z,w]\mid z_iw_i = 0 \text{ for all }i\big\}$$
is called the {\bf extended core} of $\MA$, and we will denote it $\mathcal{L}_{ext}(\A)$.
The extended core is $\Td$-equivariantly isomorphic to a union of affine toric varieties,
with moment polytopes equal to the closures of the components of $\xa$ \cite[\S 2]{HP}.
An element $[z,w]\in\mathcal{L}_{ext}(\A)$
lies on the toric divisor corresponding to the hyperplane
$H_i$ if and only if $z_i=w_i=0$ \cite[\S 3.1]{BD},
hence the generic points on $\mathcal{L}_{ext}(\A)$
consist precisely of the free $\Td$ orbits.
There is one such orbit for every component of $\xa_\R$,
hence the number of generic $\Fq$ points in the fiber $\phc^{-1}(0)$ 
is equal to $(q-1)^d \cdot r(\A)$.

Fix a flat $F$, and let $\Trkf$ be the image of $\TF$
in $\Td$, as in Section \ref{strat}.  
Then $\Trkf$ acts on $\MAf$, and the projection $p_F:\MA\to\MAf$
is $\Trkf$-equivariant.
Choose a point $x\in M\!(\AF)\subs\tdd$.
Recall from the proof of Lemma \ref{etneigh} that we have a map $p_F$ from an open
subset of $\MA$ to $\MAf$.  This open subset includes $\phc^{-1}(x)$, and
the restriction of $p_F$ maps $\phc^{-1}(x)$
surjectively onto $\mathcal{L}_{ext}(\Af)$.  The generic points in
$\phc^{-1}(x)$ are precisely those points which map to generic points of $\MAf$.
Choose a subset $S\subs\otn$ of size $\crk F$ such that $\{a_i\mid i\in F\cup S\}$
is a spanning set for $\td$, and let $T'$ be the image in $\Td$ of the coordinate torus $T^S$.
Then $T'$ acts freely on the generic fibers of $p_F$, and by dimension count, this
action is transitive as well.  Hence
$$\nu_{\phc^{-1}(x)}(q) = \nu_{T'}(q)\cdot \nu_{\Trkf}(q)\cdot r(\Af) = (q-1)^d\cdot r(\Af).$$
Summing over all $x\in M\!(\AF)$ contributes a factor of $\chi_{\AF}(q)$,
and summing over all flats $F$ yields the desired formula.
\end{proof}

\vspace{\baselineskip}
Since every stratum of $\MA$ is itself the generic stratum of some hypertoric
variety, Proposition \ref{count} implies that $\MA$ has polynomial point count,
and furthermore gives us a combinatorial formula
for counting points on each stratum.
Let $\scrT$ be the class of all hypertoric varieties corresponding
to central arrangements.  By Lemma \ref{etneigh}, $\scrT$ satisfies
condition (1).  To see that the intersection cohomology groups
of hypertoric varieties are pure, we use the decomposition theorem
of \cite[6.2.5]{BBD}, which will be discussed further in Section \ref{dc}.
This theorem implies that
the projective map
of Corollary \ref{semismall} induces an injection of $\IH(\MA)$
into the $\ell$-adic \'etale cohomology group of $\MAt_{\Fpbar}$,
which we observed was pure in Section \ref{smoothsec}.
This injection is equivariant with respect to the Frobenius action,
therefore $\IH(\MA)$ is pure as well.
Let $P_{\A}(q) = P_{\MA}(q)$.
Combining Theorem \ref{affine} with Proposition \ref{count},
Lemma \ref{etneigh}, and the
isomorphism \eqref{glo-loc-iso}
produces the following equation:
\begin{equation}\label{klrecursion}
q^{2d}\cdot P_{\A}(q^{-1}) = \sum_F P_{\Af}(q)\cdot (q-1)^{\crk F}
\cdot\sum_{G\supseteq F} \chi_{\A^G}(q)\cdot r(\AF_G),
\end{equation}
where $\AF_G = (\AF)_G = (\A_G)^F$.

\begin{theorem}\label{ip}
The intersection cohomology Poincar\'e polynomial of $\MA$ coincides
with the $h$-polynomial of $\bcsda$, that is $\paq = \hsaq$.
\end{theorem}

\begin{proof}
The polynomial $\paq$ is completely determined by
Equation \eqref{klrecursion} and the fact that $\deg\paq\leq d-1$.
It therefore suffices to prove the recursion
$$q^{2d} \cdot h_{\!\A}^{br}(q^{-1}) = \sum_F h^{br}_{\!\Af}(q)\cdot (q-1)^{\crk F}
\cdot\sum_{G\supseteq F} \chi_{\A^G}(q)\cdot r(\AF_G).$$
We proceed by expressing every piece of the equation in terms of the M\"obius
function\footnote{The M\"obius function
should not be confused with the moment map, for which we have also used the symbol $\mu$.
The M\"obius function will not appear in this paper outside of the proof of Theorem \ref{ip}.}
$$\mu:L(\A)\times L(\A)\to\Z.$$
The function $\mu$ is defined by the recursion
$$\mu(F,G)=0\text{ unless $F\subs G,\hspace{10pt}$ and if $F\subs G$, then }
\hspace{5pt}\sum_{F\subs H\subs G}\mu(H,G) = \delta(F,G),$$
where $\delta$ is the Kronecker delta function.
Let $\mu(F) = \mu(\emptyset,F)$ for all flats $F\in L(\A)$.
We may express all relevant polynomials in terms of the M\"obius
function as follows \cite[\S 7.4]{Bj}\footnote{Note that Bj\"orner's
definition of the $h$-polynomial differs from ours in that the
order of the coefficients is reversed.}:
$$\chi_{\A}(q) = \sum_F\mu(F)q^{\crk F},\hspace{10pt}
r(\A) = (-1)^{\rk\A}\chi_{\A}(-1)
= \sum_F(-1)^{\rk F}\mu(F),$$
$$\text{and}\hspace{10pt}\hsaq = (-q)^{\rk\A}\chi_{\A}(1-q^{-1})
= (-1)^{\rk\A}\sum_F\mu(F)q^{\rk F}(q-1)^{\crk F}.$$
It follows that
\begin{eqnarray*}
&& \sum_F h^{br}_{\!\Af}(q)\cdot (q-1)^{\crk F}
\cdot\sum_{G\supseteq F} \chi_{\A^G}(q)\cdot r(\AF_G)\\
&\!\!\!\!\!\!=&\!\!\!\!\!\!\!\!\!\sum_{H\subs F\subs J \subs G\subs I}
\!\!\!\!\!\!(-1)^{\rk F}\mu(H)q^{\rk H}(q-1)^{\rk F-\rk H}
\cdot (q-1)^{\crk F}
\cdot \mu(G,I)q^{\crk I}
\cdot \mu(F,J)(-1)^{\rk J - \rk F}\\
&\!\!\!\!\!\!=&\!\!\!\!\!\!\!\!\!\sum_{H\subs F\subs J \subs G\subs I}
\!\!\!\!\!\!\mu(H)q^{\rk H}(q-1)^{\crk H}
\cdot \mu(G,I)q^{\crk I}
\cdot \mu(F,J)(-1)^{\rk J}.
\end{eqnarray*}
We now apply the recursive definition of $\mu$ twice,
once to the sum over $F$ and once to the sum over $G$,
to obtain a sum over a single variable.  This yields
\begin{eqnarray*}
\sum_F(-1)^{\rk F}\mu(F)q^{\rk F+\crk F}(q-1)^{\crk F}
&=& q^{\rk\A}\sum_F(-1)^{\rk F}\mu(F)(q-1)^{\crk F}\\
&=& (-1)^{\rk\A} q^{2\rk\A}\sum_F\mu(F)q^{-\rk F}(q^{-1}-1)^{\crk F}\\
&=& q^{2\rk\A}\cdot h_{\!\A}^{br}(q^{-1}).
\end{eqnarray*}
Since $d=\rk\A$, this completes the recursion, and therefore the proof of Theorem \ref{ip}.
\end{proof}
\end{section}

\begin{section}{The KRS convolution formula}\label{dc}
In this section we use the decomposition theorem of \cite[6.2.5]{BBD} to
compare the intersection cohomology groups of $\MA$ to the ordinary
cohomology groups of its resolution $\MAt$. By the results of Sections
\ref{smoothsec} and \ref{cp}, we know that the formula that we obtain will
involve the $h$-numbers of matroid complexes and their broken circuit
complexes. In fact, this formula turns out to be a special case of the
Kook-Reiner-Stanton convolution formula, which is proven by combinatorial means in 
\cite{KRS}.\footnote{We thank Ed Swartz for this observation.}

Rather than stating the decomposition theorem for arbitrary
projective maps $f:X\to Y$, we specialize to the case where $X$ is
a complex orbifold, $Y = \sqcup_{\b\in B}\yb$ is a stratified
complex variety, and $f$ is semismall.  For the remainder of the paper
we will always work over the complex numbers, and omit the subscript $\C$.
Let $n_{\b}$ be the
codimension on $\yb$ inside of $Y$.

\begin{proposition}\label{cohomdecom}{\em\cite[4]{BM}, \cite[5.4]{Gi}}
There is a direct sum decomposition
$$H^*(X)=\bigoplus_{\b\in B}I\! H^*(\overline\yb;\eb)[-n_\b],$$
where $\eb$ is the local system
$R^{n_\b}f_*\C_{X_\b}$.
If this local system is trivial for all $\b$, then
$$H^*(X)=\bigoplus_{\b\in B}\IH(\overline{Y_\b})\otimes H^{n_{\!\b}}(\pi^{-1}(y_\b)),$$
where $y_\b\in Y_\b$, and $H^{n_{\!\b}}(\pi^{-1}(Y_\b))$ 
is understood to lie in degree $n_\b$.
\end{proposition}

We begin by showing that in the hypertoric setting, the
affinization map induces trivial local systems.

\begin{proposition}\label{trivsys}
For any flat $F$ of $\A$, the local system $\ef$ on
$\MAFopen$ induced by the affinization map $\pi:\MAt\to\MA$ is
trivial.
\end{proposition}

\begin{proof}
The stratum $\MAFopen$ admits a free action of $\Tcrkf=\Td/\Trkf$,
where $\Trkf$ is the image of the coordinate torus $\TF\subs\Tn$,
as in Section \ref{strat}. Let $\Tcrkf_\R$ be the compact real form
of $\Tcrkf$. The local system $\ef$ is naturally
$\Tcrkf_\R$-equivariant, which means that it may be pulled back
from a local system on the quotient. The quotient space
$\MAFopen/\Tcrkf_\R$ is homeomorphic, via the {\em hyperk\"ahler}
moment map, to the space $$H_F^{\R}\oplus
H_F^{\C}\smallsetminus\bigcup_{i\in F^c}H_{F\cup\{i\}}^\R \oplus
H_{F\cup\{i\}}^\C\hspace{10pt}\cite[\S 3.1]{BD}.$$ Since we are
removing linear subspaces of real codimension three, the resulting
space is simply connected, thus all of its local systems are
trivial.
\end{proof}

\vspace{\baselineskip}
By Corollary \ref{loctriv}, Corollary \ref{semismall}, Theorem \ref{cohomdecom},
and Proposition \ref{trivsys}, we obtain the following isomorphism:
\begin{equation}\label{htdecom}
H^*\big(\MAt\big)\cong \displaystyle{\bigoplus_{F\in L(\A)}
\IH\big(\MAF\big)\otimes H^{2\rk F}\!\big(\LAft\big).}
\end{equation}

\begin{remark}
The core $\LAft$ has a component of dimension $\rk F$ if and only if
$F$ is coloop-free (Remark \ref{sl}), hence only coloop-free
flats contribute to the right hand side.
\end{remark}

By Theorems \ref{sp} and \ref{ip}, Equation \eqref{htdecom} translates 
into the following recursion.

\begin{corollary}\label{poindecom}
$\displaystyle{\haq=\sum_{F\in L(\A)}h_{\!\AF}^{br}(q)\cdot h_{\rk
F}(\Af)\, q^{\rk F}}$.
\end{corollary}

\begin{remark}\label{tutte}
The {\bf Tutte polynomial} $T_\A(x,y)$ of $\da$ is
a bivariate polynomial invariant of matroids with
several combinatorially significant specializations; 
see, for example, {\em\cite[7.12~$\&$~7.15]{Bj}}.  
In particular, we have
$$\haq = q^{\rk\A}T_\A(q^{-1},1),$$
$$\hsaq=q^{\rk\A}T_\A(q^{-1},0),$$
$$\text{and}\hspace{10pt}h_{\rk\A}(\Delta_\A)=T_{\A}(0,1).$$
Corollary \ref{poindecom} is the specialization at $x = q^{-1}$ and $y=1$
of the following recursion {\em\cite[1]{KRS}:}
$$T_\A(x,y)=\sum_{F\in L(\A)}T_{\AF}(x,0)\cdot T_{\Af}(0,y).$$
\end{remark}
\end{section}

\begin{section}{Cohomology rings}\label{cr}
In Sections \ref{smoothsec} and \ref{cp}, we computed the Betti numbers of
$\MAt$ and $\MA$.  In this section, we discuss the equivariant cohomology ring
$H^*_{\!\Td}\!\big(\MAt;\C\big)$, and the equivariant intersection cohomology
group $\IHT\big(\MA;\C\big)$.  (As in Section \ref{dc}, we will consider all of our varieties
exclusively over the complex numbers.)
In general, while ordinary cohomology is a ring, intersection cohomology groups
have no naturally defined ring structure.
Nonetheless, in the case of a hypertoric variety defined by a unimodular arrangement, we conjecture that
that there is a natural isomorphism between $\IHT\big(\MA;\C\big)$
and a combinatorially defined ring $R(\A)$,
and furthermore that the multiplication in this ring has a geometric interpretation.

We begin with the ring $H^*_{\!\Td}\!\big(\MAt;\C\big)$, which was computed
independently in \cite{Ko} and \cite{HS}.
Hausel and Sturmfels observed that this ring
is isomorphic to the Stanley-Reisner ring of the matroid complex $\da$.

\begin{theorem}\label{smooth}{\em\cite[5.1]{HS},\cite[3.1]{Ko}}
There are natural ring isomorphisms
$H^*_{\!\Td}\!\big(\MAt;\C\big)\cong\SR(\DA)$
and $H^*\!\big(\MAt;\C\big)\cong\SR_0(\DA)$
that reduce degrees by half.
\end{theorem}

We now proceed to the singular case.  
The intersection cohomology group $\IHT\big(\MA;\C\big)$
is a module over $H^*_{\!\Td}(pt)$, and, as in the smooth case,
it is a free module \cite[14.1(1)]{GKM}.
We {\em could} therefore identify $\IHT\big(\MA;\C\big)$ with
the Stanley-Reisner ring of the broken circuit complex $\bcsda$, since this ring
is also a free module over a polynomial ring of dimension $d$, and has
the correct Hilbert series by Theorem \ref{ip}.  
We submit, however, that this identification
would not be natural.  One immediate objection is that the broken circuit complex
depends on a choice of ordering of the set $\otn$, while the
hypertoric variety $\MA$ does not. Instead we introduce the
the following ring,
which (if $\A$ is unimodular) is a deformation of the Stanley-Reisner ring
of the broken circuit complex for {\em any} choice of ordering.

For any circuit $C$ of $\A$, 
there exist nonzero integers $\{\lambda_i\mid i\in C\}$,
unique up to simultaneous scaling, such that $\sum_{i\in {C}}\lambda_i a_i = 0$.
Let 
\begin{equation}\label{fc}
\fc := \sum_{i\in {C}}
\,\,\operatorname{sign}(\lambda_i)\hs\!\cdot \!\!\!\prod_{j\in{C}\setminus\{i\}}\!\! e_j,
\end{equation}
where $\operatorname{sign}(\lambda) = \frac{\lambda}{|\lambda|}$.
Note that the leading term of $\fc$ with respect to an ordering $\sigma$
is equal to a multiple of the $\sigma$-broken circuit monomial corresponding
to the circuit $C$.
We now define the rings
$$R(\A) := \C[e_1,\ldots,e_n]\Big/\Big\langle \,\fc\bigmid C\text{ a circuit}\,\Big\rangle
\hspace{10pt}\text{and}\hspace{10pt}
R_0(\A) := R(\A)\otimes_{\Sym\tdd}\C,$$
where $\Sym\tdd$ acts on $R(\A)$ via the inclusion $\tdd\hookto\tnd = \C\{e_1,\ldots,e_n\}$.
If $\A$ is unimodular, then $\lambda_i$ may be chosen to be plus or minus $1$ for all $i$.
In this case, $\operatorname{sign}(\lambda_i) = \lambda_i$, and
the ring $R(\A)$ may be interpreted as the subring of rational functions
on $\tdd$ generated by the inverses of the linear functions $\{a_i\}$ \cite[\S 1]{PS}.

\begin{theorem}\label{ugb}{\em\cite[4,7]{PS}}
If $\A$ is unimodular, then
the set $\{\,\fc\mid{C}\text{ a circuit}\,\}$ 
is a universal Gr\"obner basis for
the ideal of relations in $R(\A)$.  
Thus for any ordering $\sigma$, $R(\A)$ is a deformation
of the Stanley-Reisner ring $\SR(\bcsda)$.
Furthermore, $R(\A)$ is a free module over $\Sym\tdd$, hence
$R_0(\A)$ has Hilbert series $\hsaq$.
\end{theorem}

Theorem \ref{ugb} fails if $\A$ is not unimodular.  If, for example, we
take $\A$ to consist of four lines in the plane, then we will find that
the Krull dimension of $R(\A)$ is $1$, and the Krull dimension of $\SR(\bcsda)$
is $2$.
Theorems \ref{ip} and \ref{ugb} together give us the following corollary.

\begin{corollary}\label{ring}
If $\A$ is unimodular, then
there exist graded vector space isomorphisms
$$\IHT\big(\MA;\C\big) \cong R(\A)\hspace{8pt}\text{and}\hspace{8pt}
\IH \big(\MA;\C\big) \cong R_0(\A)$$
that reduce degrees by half.
\end{corollary}

\begin{conjecture}\label{nat}
These isomorphisms are natural, 
and the multiplication in the ring $R(\A)$ 
may be interpreted as an intersection pairing on $\MA$.
\end{conjecture}

We conclude this section by providing evidence for Conjecture \ref{nat}.
Although intersection cohomology is not functorial with respect to
arbitrary maps, any map between stratified spaces with the property that
perverse $i$-chains on the source push forward to perverse $i$-chains on
the target induces a pullback in intersection cohomology.  In \cite[\S 5.4]{GM2},
the authors define the notion of a {\bf normally nonsingular map}, and prove that
such maps have this property.

For any flat $F$ of $\A$, consider the map
$s_F:\MAf\to\MA$ which exhibits $\MAf$ as a normal slice to the stratum
$\MAFopen\subs\MA$.  The stratification of $\MAf$ is pulled back from the
stratification of $\MA$, hence $s_F$ is a normally nonsingluar inclusion.
It is also $\Td$-equivariant, where $\Td$ acts on
$\MAf$ by first projecting onto $\Trkf$.
This means that if Conjecture \ref{nat} is correct, then
we should expect to find a map of rings
from $R(\A)$ to $R(\Af)$.

The ring $R(\A)$ is defined to be a quotient of
$\C[e_1,\ldots,e_n]$, while $R(\Af)$ is a quotient of
$\C[e_i]_{i\in F}$. Let $s(e_i) = e_i$ for all $i\in F$, and
and zero otherwise. To check that
$s$ is well defined, we must examine its behavior on the
element $\fc$ for every circuit ${C}$ of $\A$.  If
${C}$ is contained in $F$, then it is also a circuit of
$\Af$, and is therefore zero in $R(\Af)$.  If ${C}$ is not
contained in $F$, then the fact that $F$ is a flat implies that
$|{C}\cap F^c|\geq 2$, and therefore that $s(\fc)=0$.
Thus there is a map from $R(\A)$ to $R(\Af)$ arising naturally
from the combinatorial perspective.
The inclusion $\MAF\hookto\MA$ of Proposition
\ref{ressub} does not push perverse chains forward to perverse chains, therefore
we do not expect to find a natural map from $R(\A)$ to $R(\AF)$, and indeed we cannot.

The inclusion of the generic stratum $\MAopen$ into $\MA$ is a good map,
therefore there is a natural restriction homomorphism 
$\rho$ from $\IHT\big(\MA;\C\big)$ 
to $H^*_{\Td}(\MAopen;\C)$.
We know that $\IHT\big(\MA;\C\big)$ is
concentrated in even degree, and by the discussion that follows,
so is $H^*_{\Td}(\MAopen;\C)$.  A spectral sequence argument using parity
vanishing along the lines of \cite[3.4.1]{BGS} shows that the relative
equivariant intersection cohomology $\IHT\big(\MA,\MAopen;\C\big)$ also
has no odd cohomology, thus the long exact sequence for intersection
cohomology tells us that $\rho$ is surjective.

The $\Td$-equivariant cohomology of $\MAopen$ is 
isomorphic to the $\Td_\R$-equivariant
cohomology, where $\Td_\R$ is the compact real form of $\Td$.
As we observed in the proof of Proposition \ref{trivsys}, 
$\Td_\R$ acts freely on $\MAopen$
with quotient $$N(\A) := \tdd_{\R}\oplus\tdd_{\C}\smallsetminus
\bigcup_{i=1}^n\,\, H_i^\R \oplus
H_i^\C\hspace{10pt}\cite[\S 3.1]{BD}.$$
This space is the complement of a collection of 
codimension 3 real linear subspaces
of the vector space $\tdd_{\R}\oplus\tdd_{\C}$, 
with intersection lattice identical
to that of $\A$.  The equivariant cohomology ring of $\MAopen$ is 
isomorphic to the ordinary cohomology ring of $N(\A)$, 
which is shown in \cite[5.6]{dLS} to be isomorphic to  
$R(\A)\big/\langle e_1^2,\ldots,e_n^2\rangle$,
with $\deg e_i = 2$ for all $i$.\footnote{The 
relations $e_i^2$ are omitted
from the presentation in \cite[5.6]{dLS}, but this is only a typo.
This ring is also studied in \cite[2.2]{Co}.}
Thus $R(\A)$ surjects naturally onto 
$H^*_{\Td}(\MAopen;\C)$, providing further evidence for
Conjecture \ref{nat}.
\end{section}

\footnotesize{

}


\begin{thebibliography}{10}

\bibitem[At]{At}
C.~Athanasiadis.
\newblock Characteristic polynomials of subspace arrangements and finite fields.
\newblock {\em Adv. Math.} 122 (1996), no. 2, 193--233.

\bibitem[BBD]{BBD}
A. A. Beilinson, J. Bernstein and P. Deligne.
\newblock {\em Faisceaux Pervers},
\newblock in {\it Analysis and topology on singular spaces, I (Luminy, 1981)},
5--171, Ast\'erisque, 100, Soc. Math. France, Paris, 1982.


\bibitem[BD]{BD}
R.~Bielawski and A.~Dancer.
\newblock The geometry and topology of toric hyperk\"ahler manifolds.
\newblock {\em Comm. Anal. Geom.} 8 (2000), 727--760.


\bibitem[BGS]{BGS}
A.~Beilinson, V.~Ginzburg, W.~Soergel.
\newblock Koszul duality patterns in representation theory.
\newblock {\em J. Amer. Math. Soc.} 9 (1996) no. 2, 473--527.

\bibitem[Bj]{Bj}
A.~Bj\"orner.
\newblock The homology and shellability of matroids and geometric lattices.
\newblock {\em Matroid applications}, 226--283,
\newblock Encyclopedia Math. Appl., 40, Cambridge Univ. Press, Cambridge, 1992.


\bibitem[BM]{BM}
W.~Borho and R.~MacPherson.
\newblock Repr\'esentations des groupes de Weyl et homologie
d'intersection pour les vari\'et\'es nilpotentes.  
\newblock C. R. Acad. Sci. Paris Sér. I Math. 292 (1981), no. 15, 707--710.

\bibitem[BLR]{BLR}
S.~Bosch, W.~L\"utkebohmert, and M.~Raynaud.
\newblock {\em N\'eron models.}
\newblock Ergebnisse der Mathematik und ihrer Grenzgebiete (3), 21.
\newblock Springer-Verlag, Berlin, 1990.

\bibitem[BP]{BP}
V.~Buchstaber and T.~Panov. 
\newblock {\em Torus actions and their applications in topology and combinatorics.} 
\newblock University Lecture Series, 24. 
\newblock American Mathematical Society, Providence, RI, 2002.


\bibitem[Co]{Co}
R.~Cordovil.
\newblock A commutative algebra for oriented matroids.
\newblock {\em Discrete Comput. Geom.} 27 (2002), 73--84.

\bibitem[CB]{CB}
W.~Crawley-Boevey.
\newblock Geometry of the moment map for representations of quivers.
\newblock {\em Compositio Math.} 126 (2001), no. 3, 257--293.

\bibitem[CBVdB]{CBVdB}
W.~Crawley-Boevey and M.~Van den Bergh.
\newblock Absolutely indecomposable representations and Kac-Moody Lie algebras.
\newblock With an appendix by Hiraku Nakajima.
\newblock {\em Invent. Math.} 155 (2004), no. 3, 537--559.


\bibitem[dLS]{dLS}
M.~de Longueville and C.A.~Schultz.
\newblock The cohomology rings of complements of subspace arrangements.
\newblock {\em Math. Ann.} 319 (2001), no. 4, 625--646.

\bibitem[Gi]{Gi}
V. Ginzburg.
\newblock {\em Geometric methods in the representation theory of Hecke
algebras}.
\newblock Notes by Vladimir Baranovsky
\newblock NATO Adv. Sci. Inst. Ser. C Math. Phys. Sci., 514,
\newblock {\em Representation theories and algebraic geometry (Montreal, PQ,
1997),} 127--183, Kluwer Acad. Publ., Dordrecht, 1998.

\bibitem[GKM]{GKM}
M.~Goresky, R.~Kottwitz, and R.~MacPherson.
\newblock Equivariant cohomology, Koszul duality, and the localization theorem.
\newblock {\em Invent. Math.} 131 (1998), no. 1, 25--83.

\bibitem[GM1]{GM1}
M.~Goresky and R.~MacPherson.
\newblock Intersection homology theory.
\newblock {\em Topology} 19 (1980), no.~2, 135--162.

\bibitem[GM2]{GM2}
M.~Goresky and R.~MacPherson.
\newblock Intersection homology. II.
\newblock {\em Invent. Math.} 72 (1983), no.~1, 77--129.


\bibitem[HP]{HP}
M.~Harada and N.~Proudfoot.
\newblock Properties of the residual circle action on a hypertoric variety.
\newblock {\em Pacific J. Math.} 214 (2004), no. 2, 263--284.

\bibitem[Ha]{Ha}
T.~Hausel.
\newblock Quaternionic geometry of matroids.
\newblock arXiv:  math.AG/0308146.

\bibitem[HS]{HS}
T.~Hausel and B.~Sturmfels.
\newblock Toric hyperk\"ahler varieties.
\newblock {\em Doc. Math.}  7  (2002), 495--534 (electronic).
\newblock arXiv:  math.AG/0203096.

\bibitem[KL]{KL}
D.~Kazhdan and G.~Lusztig.
\newblock Schubert varieties and Poincar\'e duality.
\newblock Geometry of the Laplace operator
(Proc. Sympos. Pure Math., Univ. Hawaii, Honolulu, Hawaii, 1979), pp. 185--203,
Proc. Sympos. Pure Math., XXXVI, Amer. Math. Soc., Providence, R.I., 1980.

\bibitem[KW]{KW}
R.~Kiehl and R.~Weissauer.
\newblock {\em Weil conjectures, perverse sheaves and
$\ell$-adic Fourier transform}.
\newblock Springer-Verlag, New York, 2001.


\bibitem[Ko]{Ko}
H.~Konno.
\newblock Cohomology rings of toric hyperk\"ahler manifolds.
\newblock {\em Int. J. of Math.} 11 (2000)
no. 8, 1001--1026.

\bibitem[KRS]{KRS}
W.~Kook, V.~Reiner and D.~Stanton.
\newblock A convolution formula for the Tutte polynomial.
\newblock {\em J. Combin. Theory} Ser. B 76 (1999), no. 2, 297--300.
\newblock arXiv:  math.CO/9712232.

\bibitem[Lu]{Lu}
G.~Lusztig.
\newblock Fermionic form and Betti numbers.
\newblock arXiv:  math.QA/0005010.

%

\bibitem[Na]{Na}
H.~Nakajima.
\newblock Quiver varieties and finite-dimensional representations of quantum affine algebras.
\newblock {\em J. Amer. Math. Soc.} 14 (2001) no. 1, 145--238.


\bibitem[P1]{P1}
N.~Proudfoot.
\newblock Hyperk\"ahler analogues of K\"ahler quotients.
\newblock Ph.D. Thesis, U.C. Berkeley, Spring 2004.
\newblock arXiv:  math.AG/0405233.

\bibitem[P2]{P2}
N.~Proudfoot.
\newblock Geometric invariant theory and projective toric varieties.
\newblock {\em Snowbird Lectures in Algebraic Geometry},
\newblock Contemp. Math. 388, {\em Amer. Math. Soc., Providence, RI}, 2005. 

\bibitem[PS]{PS}
N.~Proudfoot and D.~Speyer.
\newblock A broken circuit ring.
\newblock arXiv:  math.CO/0410069.

\bibitem[SL]{SL}
R.~Sjamaar and E.~Lerman.
\newblock Stratified symplectic spaces and reduction.
\newblock {\em Ann. of Math.} (2) 134 (1991), no. 2, 375--422.

\bibitem[S1]{S1}
R.~Stanley.
\newblock The number of faces of a simplicial convex polytope.
\newblock {\em Adv. in Math.} 35 (1980), no. 3, 236--238.

\bibitem[S2]{S2}
R.~Stanley.
\newblock Generalized $H$-vectors, intersection cohomology of toric varieties, and related results.
\newblock {\em Commutative algebra and combinatorics}
\newblock (Kyoto, 1985), 187--213, Adv. Stud. Pure Math., 11, North-Holland, Amsterdam, 1987.
 
\bibitem[S3]{S3}
R.~Stanley.
\newblock {\em Combinatorics and commutative algebra}.
\newblock Progress in Mathematics, 41.
\newblock Birkhäuser Boston, Inc., Boston, MA, 1983.

\end{thebibliography}
\end{document}